\documentclass[12pt,a4paper,reqno]{amsart}
\usepackage{graphics,epic}
\usepackage{amsmath,amssymb, amsthm}
\usepackage[all,2cell]{xy}

\newtheorem{theorem}{Theorem}[section]
\newtheorem*{theorem*}{Theorem}

\newtheorem{lemma}[theorem]{Lemma}
\newtheorem{proposition}[theorem]{Proposition}

\newtheorem*{conjecture*}{Conjecture}

\newtheorem{remark}[theorem]{Remark}

\newcommand{\ca}{{\mathcal A}}

\newcommand{\cc}{{\mathcal C}}
\newcommand{\cd}{{\mathcal D}}

\newcommand{\cf}{{\mathcal F}}

\newcommand{\ct}{{\mathcal T}}

\newcommand{\cz}{{\mathcal Z}}

\newcommand{\bt}{\boxtimes}

\renewcommand{\hat}[1]{\widehat{#1}}

\newcommand{\ot}{\otimes}

\newcommand{\id}{{\rm id}}
\newcommand{\Hom}{{\rm Hom}\,}

\newcommand{\End}{{\rm End}\,}

\newcommand{\FPdim}{{\rm FPdim}\,}

\newcommand{\Z}{\mathbb{Z}}

\newcommand{\C}{\mathbb{C}}

\def\wt{{\rm wt}}

\def\C{{\mathbb C}}

\def\Z{{\mathbb Z}}

\def\Y{{\mathcal Y}}

\def\1{{\bf 1}}

\def \End{{\rm End}}
\def \Hom{{\rm Hom}}

\def \pf{\noindent {\bf Proof: \,}}
\def\theequation{5.\arabic{equation}}

\def \l{\lambda}
\def \w{\omega}

\begin{document}

\title[Mirror extensions of rational vertex operator algebras]{Mirror extensions of rational vertex operator algebras}

\author{Xingjun Lin}
\address{Xingjun Lin, Institute of Mathematics, Academia Sinica, Taipei 10617, Taiwan}
\email{linxingjun@math.sinica.edu.tw}

\begin{abstract}
In this paper,  mirror extensions of rational vertex operator algebras are considered. The mirror extension conjecture is proved.
\end{abstract}
\maketitle
\section{Introduction \label{intro}}
\def\theequation{1.\arabic{equation}}
\setcounter{equation}{0}
Mirror extensions were first studied in the context of conformal nets, new conformal nets were obtained and it was further proved that these conformal nets can not be obtained by cosets, orbifolds and simple current extensions  \cite{X1}. The mirror extensions of vertex operator algebras were  studied in \cite{DJX},  mirror extensions of specific vertex operator algebras were obtained and it was  conjectured that the mirror extensions exist for general vertex operator algebras satisfying some conditions. More precisely,   let $(U, Y, \1, \w)$ be a vertex operator algebra and
$(V, Y, \1, \w')$ be a vertex operator subalgebra of $U$ such
that $\w'\in U_2$ and $\w_2\w'=0$, denote $(V^c, Y, \1, \w-\w')$ the commutant vertex operator algebra of $V$ in $U$ (see \cite{FZ}).
 Assume further that $(V, Y, \1, \w)$, $(V^c, Y, \1, \w-\w')$,
$(U, Y, \1, \w')$ are rational and $C_2$-cofinite vertex operator algebras such that
  $(V^c)^c=V$ and $U$ has the following decomposition as  $V\otimes V^c$-module: $$U=V\ot V^c\oplus (\oplus_{i=1}^nM^i\ot N^i),$$ where $M^i$'s and $N^i$'s are irreducible modules for $V$ and $V^c$ respectively. Let $$V^e=V\oplus (\oplus_{i=1}^nm_iM^i)$$ and $$(V^c)^e=V^c\oplus (\oplus_{i=1}^nm_i(N^i)'),$$ where $m_i$'s are nonnegative integers. It was conjectured in \cite{DJX} that if there is a vertex operator algebra structure on $V^e$ such that $V^e$ is an extension  vertex operator algebra of $V$, then   $(V^c)^e$ has a vertex operator algebra structure such that $(V^c)^e$ is an extension vertex operator algebra of $V^c$, the so-called  mirror extension vertex operator algebra of $V^c$. The main goal of this paper is to prove this conjecture.

In our proof of the mirror extension conjecture, we in fact prove that if $(V, Y, \1, \w)$, $(V^c, Y, \1, \w-\w')$,
$(U, Y, \1, \w')$ are rational and $C_2$-cofinite vertex operator algebras such that
  $(V^c)^c=V$,  then $U$ has the following decomposition as  $V\otimes V^c$-module:
  \begin{align}\label{ivoa}
  &U=V\ot V^c\oplus (\oplus_{i=1}^nM^i\ot N^i)
   \end{align}
   where  $M^i$'s and $N^i$'s are inequivalent irreducible modules for $V$ and $V^c$ respectively. This fact plays an important role in our proof of the mirror extension conjecture. Explicitly, let  $\cc_{V}$, $\cc_{V^c}$ be the module categories of $V$, $V^c$, respectively. We prove that
 $\{V, M^i| 1\leq i\leq n \}$ $($resp. $\{V^c,N^i| 1\leq i\leq n \})$ are closed under the tensor product of $\cc_{V}$ (resp. $\cc_{V^c}$). In particular, $\{V, M^i| 1\leq i\leq n \}$ $($resp. $\{V^c, N^i| 1\leq i\leq n \})$ forms a braided monoidal subcategory $\cc_{V}^0$ $($resp. $\cc_{V^c}^0)$ of  $\cc_{V}$ $($resp. $\cc_{V^c})$. We further prove that the decomposition (\ref{ivoa}) induces a braid-reversing equivalence between $\cc_{V}^0$ and $\cc_{V^c}^0$ such that $M^i$ maps to $(N^i)'$ for any $1\leq i\leq n$. Then the mirror extension conjecture could be proved immediately by using results in \cite{HKL}. These generalize the results in \cite{OS}, where the  level-rank duality \cite{F} was studied via tensor category and the category of integrable highest weight representations of affine lie algebras $\hat{sl}_m$ at level $n$ was proved to be  partly braid-reversing equivalent to the category of integrable highest weight representations of affine lie algebras $\hat{sl}_n$ at level $m$. From this point of view,  $V$ is dual to $V^c$.

This paper is organized as follows: In Section 2, we recall some basic facts about tensor categories and vertex operator algebras. In
Section 3, we prove our main results on mirror extensions of vertex operator algebras.

\section{Preliminaries }
\def\theequation{2.\arabic{equation}}
\setcounter{equation}{0}
\subsection{Preliminaries on tensor categories}
In this subsection, we recall some notions and facts in the theory of tensor categories from \cite{BK}, \cite{DMNO}, \cite{KO} and \cite{O}.
\subsubsection{Braided tensor categories}
Let $\cc$ 
be a monoidal category defined as in \cite{BK}  and  denote the unit object in $\cc$ by $\1_{\cc}$.
 For  an object $U$ in  $\mathcal{C}$, 
 a {\em (right) dual} of $U$ is an object $U^\ast$ with two morphisms
$$e_{U}: U^{\ast}\otimes U\longrightarrow \1_{\cc},$$
$$i_{U}: \1_{\cc} \longrightarrow  U\otimes U^{\ast},$$ such that the composition $$U\stackrel{i_U\otimes \id_U}{\longrightarrow} U\otimes U^{\ast}\otimes U\stackrel{\id_U\otimes e_U}{\longrightarrow} U$$ is equal to $\id_{U}$ and the composition $$U^{\ast}\stackrel{\id_{U^{\ast}}\otimes i_U}{\longrightarrow} U^{\ast}\otimes U\otimes U^{\ast}\stackrel{e_U\otimes \id_{U^{\ast}}}{\longrightarrow} U^{\ast}$$ is equal to $\id_{U^{\ast}}$.
Similarly, for any object $U\in \cc$ we could define the {\em left dual} $\,^*U$ of $U$.
If the dual of $U$ exist, then   the dual of $U$ is unique (see \cite{BK}). Furthermore, we have the following results which will be used frequently in this paper (see \cite{BK}).
\begin{proposition}\label{kc3}
Let $\cc$ 
 be a monoidal category and $V$ be an object in $\cc$. If $V$ has a dual $V^*$, then
\begin{align*}
&\Hom(U\ot V, W)=\Hom(U, W\ot V^*),\\
&\Hom(U, V\ot W)=\Hom(V^*\ot U,  W).
\end{align*}
\end{proposition}
Recall \cite{BK} that a monoidal category $\mathcal{C}$
is called {\em a braided tensor category} if there is a natural bifunctor isomorphism $c_{X, Y}: X\otimes Y\to Y\otimes X$ called
a {\em braiding} subject to the hexagon axioms. For every braided tensor category $(\mathcal C, c)$,
there is a {\em reversed braiding} on $\mathcal C$ given by $c_{X, Y}^{rev}=c_{Y, X}^{-1}$, a braided tensor category $\mathcal C$
endowed with the reversed braiding will be denoted $\mathcal C^{rev}$.
 Two
objects $X, Y$ in a braided tensor category $\cc$ are said to {\em mutually centralize} each other, in
the sense of \cite{M}, if $c_{X,Y} c_{Y,X}= \id_{X\otimes Y} $. For a  subcategory $\mathcal D$ of a braided tensor category
$\mathcal C$, its {\em centralizer} $\mathcal D'$
 is the full subcategory of $\mathcal C$ consisting of all objects that centralize
every object of $\mathcal D$.

For a monoidal category $\mathcal C$, the {\em Drinfeld center} $\cz(\mathcal C)$ is defined to be the
category whose objects are pairs $(V, \gamma_{-,V})$ for $V$ an object of $\mathcal C$, and $\gamma_{-,V}$
is a natural family of isomorphisms  $\gamma_{X, V}: X \otimes V\to   V \otimes X$ defined for all objects $X\in \mathcal C$ such that for all objects $X, Y\in \cc$,$$\gamma_{X\ot Y, V}=(\gamma_{X, V}\ot \id_Y)(\id_X\ot \gamma_{Y, V}).$$A morphism from $(V, \gamma_{-,V})$ to $(W, \gamma_{-,W})$ is a morphism $f:V\to W$ in $\cc$ such that for every object $X\in \cc$, $$(f\ot \id_X)\gamma_{X, V}=\gamma_{X, W}(\id_{X}\ot f).$$
The following theorem was proved in \cite{Ka}.
\begin{theorem}\label{kvoa10}
For a monoidal category $\mathcal C$,  $\cz(\cc)$ is  a braided tensor category and the braiding is given by $$\gamma_{V, W}: (V, \gamma_{-, V})\to (W, \gamma_{-,W}).$$
\end{theorem}

 Let $\cc,\cd$ be monoidal categories. If $\cc,\cd$ are braided tensor categories, then a monoidal functor $\mathcal T : \mathcal C\to \mathcal D $ is called {\em braided} if  $\ct$ preserves the
braiding; A braided monoidal functor $\mathcal T : \mathcal C\to \mathcal D $ is said
to be a {\em braided equivalence of categories} if it is  an equivalence of the underlying categories (see  \cite{BK}).
A {\em braid-reversing equivalence} $ \mathcal T :\mathcal C\to \mathcal D $ is a braided equivalence $\mathcal C\to \mathcal D^{rev}$.
Given a monoidal functor $\mathcal T : \mathcal C\to \mathcal D$ with $\mathcal C$ braided, we say $\mathcal T$ is {\em central} if it
factors through the forgetful functor $\cz(\mathcal D)\to \mathcal  D$ via a braided monoidal functor $\mathcal  C \to \cz(\mathcal D)$.

If $\cc$ is a braided tensor category, we have the functors
\begin{align*}
&i : \mathcal C \to \cz(\mathcal C),\ \ \ \ X\mapsto (X, c_{-, X} ),\\
&j : \mathcal C \to  \cz(\mathcal C),\ \ \ \  X\mapsto (X, c_{ X, -}^{-1})
\end{align*}
are faithfully braided and braid-reversing monoidal
functors, respectively (see \cite{Ka}, \cite{M2}).

Recall \cite{BK} that a monoidal category $\cc$ is called {\em rigid} if every object in $\cc$ has right and left duals.
A   monoidal category $\cc$ is called a {\em fusion category} if $\cc$
is a semisimple  rigid monoidal category with finite dimensional spaces of morphisms,
finitely many irreducible objects and an irreducible unit object. The following result was obtained in \cite{M2}.
\begin{theorem}\label{kc5}
Let $\mathcal C$ be a braided fusion category, $i, j$ be the functors defined above. Then we have
$$i(\mathcal C)'= j(\mathcal C).$$
\end{theorem}

\subsubsection{Algebras in monoidal categories}
Let $\cc$ be a monodial category and  denote the unit object in $\cc$ by $\1_{\cc}$. An {\em algebra} in $\mathcal{C}$ is an object $A\in \mathcal{C}$ along with morphisms $\mu: A\otimes A\to A$ and $\iota_A:\1_{\cc}\hookrightarrow A$ such that the following conditions hold:

(i) {\em Associativity}. Compositions $\mu\circ (\mu\otimes \id)\circ a,\ \mu\circ (\id\otimes \mu):A\otimes (A\otimes A)\to A$ are equal, where $a$ denotes the associativity isomorphism $a: A\otimes (A\otimes A)\to (A\otimes A)\otimes A$;


(ii) {\em Unit}. Composition $\mu \circ (\iota_A\otimes A):A=\1_{\cc}\otimes A \to A$ is equal to $\id_A$;

This completes the definition of the  algebra. We will denote the  algebra just defined by $(A, \mu, \iota_A)$ or briefly by $A$.

 For  an  algebra $A$  in $\mathcal{C}$, the {\em left module category} ${}_A\cc$  is defined to be the category whose objects  are pairs $(M,\mu_M)$, where $M\in \mathcal{C}$ and $\mu_M: A\otimes M\to M$ is a morphism in $\mathcal{C}$ such that:
\begin{align*}
& \mu_M\circ (\mu\otimes \id)\circ a=\mu_M\circ (\id \otimes \mu_M):A\otimes (A\otimes  M)\to M;\\
& \mu_M(\iota_A\otimes \id)=\id: \1\otimes M\to M.
\end{align*}
The morphisms are defined by
\begin{align*}
&\Hom_{{}_A\cc}((M_1,\mu_{M_1}), (M_2,\mu_{M_2}))\\
&\ \ \ \ \ \ \ \ \ \ \ =\{\phi\in \Hom_{\mathcal{C}}(M_1, M_2)|\mu_{M_2}\circ(\id\otimes \phi)=\phi\circ\mu_{M_1}:A\otimes M_1\to M_2\}.
\end{align*}
Similarly, we could define the {\em right module category} $\cc_A$.  An  algebra $A$ in $\mathcal{C}$ is called {\em semisimple} if the category ${}_A\cc$ is semisimple. The following theorem was proved in \cite{DMNO}.
\begin{proposition}\label{kc1}
Let $\cc$ be a fusion category and $A$ be an  algebra in $\cc$. Then the following conditions are equivalent:\\
(i) The category $\cc_A$ is semisimple;\\
(ii) The category ${}_A\cc$ is semisimple.
\end{proposition}

 We now assume that  $\cc$ is a modular tensor category defined as in \cite{BK}. In particular, $\cc$ is a braided fusion category. An {\em etale algebra} $A$ in  $\mathcal C$ is defined to be an algebra such that the category $\cc_{A}$ is  semisimple and $A$ is {\em commutative}, that is, $\mu \circ c_{A,A}:A\otimes A\to A$ is equal to $\mu$. In particular, from Proposition \ref{kc1} we have $A$ is  semisimple, i.e., the category $_A\cc$ is semisimple.  An etale algebra $A$ is called {\em connected} if the unit object appears in $A$ with multiplicity 1.
 For a connected
 etale algebra $A\in  \mathcal C$, it was proved  that the category $\mathcal C_A$ is a fusion category
with operation $\otimes_A$ (see \cite{DMNO}). 
Let $\cf_A$ be the  {\em free module
functor} $\cf_A : \mathcal C\to \cc_A$ defined by $\cf_A(X) = X\otimes A$  (see \cite{KO}). It was proved in \cite{KO} that $\cf_A$ is a monoidal functor. Furthermore, it was proved  that $\cf_A$ has the structure of a central functor
given by $X\mapsto  (X \otimes A, c_{-,X})$ (see \cite{DMNO}).

  Recall from \cite{BK} that for any modular tensor category $\cc$ and  $U$ in $\cc$ there is a   functorial isomorphism
$$\theta_U: U\longrightarrow U$$  such that:
\begin{align*}
&\theta_{U\otimes W}=c_{W,U}c_{U,W}(\theta_{U}\otimes \theta_{W});\\
&\theta_{1}=\id;\\
&\theta_{U^\ast}=(\theta_U)^\ast,
\end{align*}
where $(\theta_U)^\ast\in \Hom(U^\ast, U^\ast)$ denotes the image of $\theta_U\in \Hom(U, U)$ under the canonical map. The following theorem was essentially proved in \cite{KO}.
\begin{theorem}\label{semisimple}
Let $\mathcal{C}$ be  a modular tensor category and $A$ be a commutative algebra in $\mathcal{C}$ such that  $\theta_A=\id$. 
If $A$ is simple as an $A$-module, then $A$ is semisimple.
\end{theorem}

In the proof of our main results, we will frequently use the following result (see \cite{KO}, \cite{O}).
\begin{theorem}\label{kc2}
 Let $\cc$ be a modular tensor category and  $A$ be an algebra in $\mathcal{C}$. For any object $U\in \cc$ and left  $A$-module $M$, we have
 $$\Hom_{A}(A\ot U, M)=\Hom(U, M),$$
 where we denote $\Hom_{_{A}\cc}$  by $\Hom_A$. Similarly, for any object $U\in \cc$ and  $A$-module $M$, we have
 $$\Hom_{\cc_A}(U\ot A, M)=\Hom(U, M).$$
\end{theorem}

\subsubsection{Frobenius-Perron dimension}
Recall from   \cite{ENO} that for a fusion category $\mathcal{C}$
there is a unique homomorphism from the Grothendieck ring of $\mathcal{C}$ to real numbers sending
each isomorphism class to a nonnegative real number. The value of this homomorphism
on the class represented by $X\in \mathcal{C}$ is called the {\em Frobenius-Perron dimension} of $X$ and
written $\FPdim_{\mathcal{C}}(X)$. For a fusion category $\cc$,  one defines its {\em Frobenius-Perron dimension} as follows:$$\FPdim \mathcal{C}=\sum_{X\in \mathcal{O}(\mathcal{C})}\FPdim_{\mathcal{C}}( X)^2,$$  where $\mathcal{O}(\mathcal{C})$ is the set of isomorphism classes of simple objects in $\mathcal{C}$. The following result was obtained in \cite{DMNO}.
\begin{proposition}\label{kc4}
Let $\cc$ a modular tensor category and  $A $ be a connected
etale algebra
in $\cc$. Then
$$(\FPdim_\cc A) \FPdim \cc_A = \FPdim \cc.$$

\end{proposition}

\subsection{Preliminaries on vertex operator algebras}
\subsubsection{Vertex operator algebras}
We first recall some basic notions and facts in the theory of vertex operator algebras from \cite{DLM1}, \cite{FHL}, \cite{FLM}, \cite{LL} and \cite{Z}.
A {\em vertex  operator algebra} is  a quadruple $(V, Y, \1, \w)$, where $V=\bigoplus_{n\in \Z}V_n$ is a $\Z$-graded $\C$-vector space such that $\dim V_n<\infty$ and $V_n=0$ for $n$ sufficiently large,  $\1\in V_{0}$ is the {\em vacuum vector} of $V$, $\w$ is the {\em conformal vector} of $V$, and $Y$ is a linear map \begin{align*}
Y:V&\to (\End V)[[x,x^{-1}]],\\
v&\mapsto Y(v,x)=\sum_{n\in\Z}v_nx^{-n-1},\,v_n\in \End V
\end{align*}
satisfying the following axioms:

(i) For any $u, v\in V$, $u_nv=0$ for $n$ sufficiently large;

(ii) $Y(\1, x)=\id_V$;

(iii) $Y(v, x)\1\in V[[x]]$ and $Y(v, x)\1|_{x=0}=v$ for any $v\in V$;

(iv) The component operators of $Y(\w,x)=\sum_{n\in \Z}L(n)x^{-n-2}$ satisfy the Virasoro algebra relation with {\em central charge} $c\in \C$:
$$[L(m), L(n)]=(m-n)L(m+n)+\frac{1}{12} (m^3-m)\delta_{m+n,0}c,$$ and $$L(0)|_{V_n}=n, n\in \Z,$$ $$\frac{d}{dx}Y(v, x)=Y(L(-1)v, x);$$

(v)  The {\em Jacobi identity}
\begin{align*}
\begin{split}
&x_{0}^{-1}\delta\left(\frac{x_{1}-x_{2}}{x_{0}}\right)Y(u,x_{1})Y(v,x_{2})-x_{0}^{-1}\delta\left(
\frac{x_{2}-x_{1}}{-x_{0}}\right)Y(v,x_{2})Y(u,x_{1})\\
&\quad=x_{2}^{-1}\delta\left(\frac{x_{1}-x_{0}}{x_{2}}\right)Y(Y(u,x_{0})v,x_{2}).
\end{split}
\end{align*}
This completes the definition of the vertex operator algebra and  we will denote the vertex operator algebra briefly by $V$. A vertex operator algebra $V$ is called of {\em CFT} type if $V=\bigoplus_{n\geq 0}V_n$ and $\dim V_0=1$.



Let $V$ be a vertex operator algebra. A {\em weak  $V$-module} $M$ is a vector space equipped
with a linear map
\begin{align*}
Y_{M}:V&\to (\End M)[[x, x^{-1}]],\\
v&\mapsto Y_{M}(v,x)=\sum_{n\in\Z}v_nx^{-n-1},\,v_n\in \End M
\end{align*}
satisfying the following conditions: For any $u\in V,\ v\in V,\ w\in M$ and $n\in \Z$,
\begin{align*}
&\ \ \ \ \ \ \ \ \ \ \ \ \ \ \ \ \ \ \ \ \ \ \ \ \ \ u_nw=0 \text{ for } n>>0;\\
&\ \ \ \ \ \ \ \ \ \ \ \ \ \ \ \ \ \ \ \ \ \ \ \ \ \ Y_M(\1, x)=\id_M;\\
\begin{split}
&x_{0}^{-1}\delta\left(\frac{x_{1}-x_{2}}{x_{0}}\right)Y_{M}(u,x_{1})Y_M(v,x_{2})-x_{0}^{-1}\delta\left(
\frac{x_{2}-x_{1}}{-x_{0}}\right)Y_M(v,x_{2})Y_M(u,x_{1})\\
&\quad=x_{2}^{-1}\delta\left(\frac{x_{1}-x_{0}}{x_{2}}\right)Y_M(Y(u,x_{0})v,x_{2}).
\end{split}
\end{align*}


A weak
 $V$-module $M$  is called an \textit{admissible $V$-module} if $M$  has a $\Z_{\geq
0}$-gradation $M=\bigoplus_{n\in\Z_{\geq 0}}M(n)$ such
that
\begin{align*}\label{AD1}
a_mM(n)\subset M(\wt{a}+n-m-1)
\end{align*}
for any homogeneous $a\in V$ and $m,\,n\in\Z$.
An admissible $V$-module $M$ is said to be
\textit{irreducible} if $M$ has no non-trivial admissible weak
$V$-submodule. When an admissible $V$-module $M$ is a
direct sum of irreducible admissible submodules, $M$ is called
\textit{completely reducible}.
A vertex operator algebra $V$ is said to be \textit{rational} if
any  admissible $V$-module is completely reducible.

A {\em  $V$-module} is a weak $V$-module $M$ which carries a $\C$-grading induced by the spectrum of $L(0)$, that is  $M=\bigoplus_{\lambda\in\C}
M_{\lambda}$ where
$M_\lambda=\{w\in M|L(0)w=\lambda w\}$. Moreover one requires that $M_\lambda$ is
finite dimensional and for fixed $\lambda\in\C$, $M_{\lambda+n}=0$
for sufficiently small integer $n$.
Similarly, a $V$-module $M$ is said to be
\textit{irreducible} if $M$ has no non-trivial
$V$-submodule and a vertex operator algebra $V$ is said to be {\em regular} if any weak $V$-module
$M$ is a direct sum of irreducible $V$-modules.

For a $V$-module $M = \bigoplus_{\l\in \mathbb{C}}{M(\l)}$, set $M'
= \bigoplus_{\lambda \in \mathbb{C}}{M_\lambda^*}$, the restricted
dual of $M$. It was proved in [FHL] that $M'$ is naturally a
$V$-module where the vertex operator denoted by $Y'$ is defined
by the property
$$\langle Y'(v, x)w', w\rangle  = \langle w', Y(e^{xL(1)}(-x^{-2})^{L(0)}v, x^{-1})w\rangle $$for $v\in V, w'\in
M'$ and $w\in M$. The $V$-module $M'$ is called the {\em contragredient
module} of $M$. A $V$-module $M$ is called {\em self dual} if $M$ is isomorphic to $M'$.

Recall that a vertex operator algebra $V$ is called {\em $C_2$-cofinite} if $\dim V/C_2(V)<\infty$, where $C_2(V)=\langle u_{-2}v|u, v\in V \rangle$. The following results were obtained in \cite{ABD}, \cite{L}.
\begin{theorem}
Let $V$ be a $CFT$ type vertex operator algebra. Then $V$ is regular if and only if $V$ is rational and $C_2$-cofinite.
\end{theorem}

 We now assume that $V$ is a  vertex operator algebra satisfying the following conditions:

(1) $V$ is simple CFT type vertex operator algebra  and is self dual;

(2) $V$ is regular.\\
Let $M$ be a weak $V$-module, we have $L(0)$ acts semisimplely on $M$, i.e. $$M=\oplus_{\lambda\in \C}M_{\lambda},$$ where $M_{\lambda}=\{w\in M| L(0)w=\lambda w\}$. If $w\in M_{\lambda}$, we write $\wt w=\lambda$ and call $w$ a {\em homogeneous vector}. In particular, if $M$ is an irreducible  $V$-module, it was proved in \cite{Z} that $$M=\oplus_{n=0}^{\infty}M_{\lambda+n},$$ for some $\lambda$ such that $M_{\lambda}\neq 0$, and $\lambda$ is defined to be  the {\em conformal weight} of $M$.
\subsubsection{Tensor product vertex operator algebra}Let $(V^1, Y, \1, \omega^1), \cdots, (V^p, Y, \1,\omega^p)$ be vertex operator algebras, the
tensor product of $(V^1, Y, \1, \omega^1), \cdots, (V^p, Y, \1,\omega^p)$ is
constructed on the tensor product vector space$$V = V^1\otimes
\cdots \otimes V^p,$$
where the vertex operator $Y$ is
defined by
$$Y(v^1\otimes \cdots \otimes v^p, x) = Y(v^1, x)\otimes
\cdots \otimes Y(v^p, x)$$for $v^i \in V^i\ (1\leq i \leq p),$
the vacuum vector is $$\1 = \1\otimes \cdots \otimes \1$$
and the conformal vector is
$$\omega = \omega^1 \otimes \cdots \otimes \1 + \cdots +\1\otimes
\cdots\otimes \omega^p.$$
Then  $(V, Y, 1, \omega)$
is a vertex operator algebra (see  \cite{FHL}, \cite{LL}).

Let $M^i$  be an admissible $V^i$-module for $i=1,...,p.$  We may construct
the tensor product admissible module $M^1\otimes \cdots \otimes M^p$ for the
tensor product vertex operator algebra $V^1\otimes \cdots \otimes
V^p$ by
$$Y_{M^1\otimes \cdots \otimes M^p}(v^1\otimes \cdots \otimes v^p,x) =
Y_{M^1}(v^1, x)\otimes \cdots \otimes Y_{M^p}(v^p, x).$$ Then  $(M^1\otimes \cdots \otimes M^p, Y_{M^1\otimes \cdots \otimes M^p})$ is an admissible $V^1\otimes
\cdots \otimes V^p$-module. Furthermore, we have the
 following results which were proved in \cite{FHL}, \cite{DMZ} and \cite{ABD}.
\begin{theorem}\label{kvoa3}
Let $V^1, \cdots, V^p$ be  vertex operator algebras. If $V^1, \cdots, V^p$ are rational, then
$V^1\otimes \cdots \otimes V^p$  is rational and any irreducible $V^1\otimes \cdots \otimes V^p$-module is a tensor
product $M^1\otimes \cdots \otimes M^p$ , where $M^1, \cdots, M^p$ are
some irreducible modules for the vertex operator algebras $V^1,
\cdots, V^p,$ respectively. If $V^1, \cdots, V^p$ are $C_2$-cofinite, then $V^1\otimes \cdots \otimes V^p$  is $C_2$-cofinite. In particular, if $V^1, \cdots, V^p$ are  vertex operator algebras satisfying conditions $(1)$ and $(2)$, then $V^1\otimes \cdots \otimes V^p$ satisfies conditions $(1)$ and $(2)$.
\end{theorem}
\subsubsection{Tensor category of $V$-modules}
First we  recall the  notions of intertwining operator and fusion
rule from [FHL]. Let $V$ be a vertex operator algebra,
 $M^1$, $M^2$, $M^3$ be weak $V$-modules. An {\em intertwining
operator} $\mathcal {Y}$ of type $\left(\begin{tabular}{c}
$M^3$\\
$M^1$ $M^2$\\
\end{tabular}\right)$ is a linear map
\begin{align*}
\mathcal
{Y}: M^1&\rightarrow \Hom(M^2, M^3)\{x\},\\
 w^1&\mapsto\mathcal {Y}(w^1, x) = \sum_{n\in \mathbb{C}}{w_n^1x^{-n-1}}
\end{align*}
satisfying the following conditions: For any $v\in V, w^1\in M^1, w^2\in M^2$ and $\lambda \in \mathbb{C},$
\begin{align*}
&  w_{n+\lambda}^1w^2 = 0 \text{ for }n>>0;\\
&  \dfrac{d}{dx}\mathcal{Y}(w^1,
x)=\mathcal{Y}(L(-1)w^1, x);\\
\begin{split}
x_0^{-1}\delta(\frac{x_1-x_2}{x_0})Y_{M^3}(v, x_1)&\mathcal
{Y}(w^1, x_2)-x_0^{-1}\delta(\frac{x_2-x_1}{-x_0})\mathcal{Y}(w^1,
x_2)Y_{M^2}(v, x_1)\\
&=x_2^{-1}\delta(\frac{x_1-x_0}{x_2})\mathcal{Y}(Y_{M^1}(v, x_0)w^1, x_2).
\end{split}
\end{align*}

Denote the vector space of intertwining operators of type $\left(\begin{tabular}{c}
$M^3$\\
$M^1$ $M^2$\\
\end{tabular}\right)$ by $\mathcal{V}_{M^1,M^2}^{M^3}$.  The dimension  of $\mathcal{V}_{M^1,M^2}^{M^3}$ is called the
{\em fusion rule} for $M^1$, $M^2$ and $M^3$, and is denoted by $N_{M^1,M^2}^{M^3}$. If $V$ is a vertex operator algebra satisfying  conditions (1) and (2), $M^1,\ M^2,\ M^3$ are $V$-modules, we have $N_{M^1,M^2}^{M^3}=N_{(M^1)',(M^2)'}^{(M^3)'}$ (see \cite{DJX1}).

For any intertwining operator $\Y$ of type  $\left(\begin{tabular}{c}
$M^3$\\
$M^1$ $M^2$\\
\end{tabular}\right)$, define the operator $\Y^*$ as follows: $$\Y^*(w^2, x)w^1=e^{xL(-1)}\Y(w^1, -x)w^2$$ for any $w^1\in M^1, w^2\in M^2$. It was proved in \cite{FHL} that $\Y^*$ is an intertwining
operators of type $\left(\begin{tabular}{c}
$M^3$\\
$M^2$ $M^1$\\
\end{tabular}\right)$. The following theorem was proved in \cite{ADL}.
\begin{theorem}\label{kvoa6}
Let $V^1, V^2$ be rational vertex operator algebras. Let $M^1 , M^2,
M^3$ be $V^1$-modules and  $N^1, N^2, N^3$ be $V^2$-modules such
that
$$\dim I_{V^1}\left(\begin{tabular}{c}
$M^3$\\
$M^1$ $M^2$\\
\end{tabular}\right)< \infty , \ \dim I_{V^2}\left(\begin{tabular}{c}
$N^3$\\
$N^1$ $N^2$\\
\end{tabular}\right)< \infty.$$
Then the linear map
$$\sigma: I_{V^1}\left(\begin{tabular}{c}
$M^3$\\
$M^1$ $M^2$\\
\end{tabular}\right)\otimes I_{V^2}\left(\begin{tabular}{c}
$N^3$\\
$N^1$ $N^2$\\
\end{tabular}\right)\rightarrow I_{V^1\otimes V^2}\left( \begin{tabular}{c}
$M^3\otimes N^3$\\
$M^1\otimes N^1$ $M^2\otimes N^2$\\
\end{tabular}\right)$$
$$\mathcal{Y}_1\otimes \mathcal{Y}_2\mapsto (\mathcal{Y}_1\otimes
\mathcal{Y}_2)(\cdot, x)$$\\is an isomorphism, where
$(\mathcal{Y}_1\otimes \mathcal{Y}_2)(\cdot, x)$ is defined by
$$(\mathcal{Y}_1\otimes \mathcal{Y}_2)(\cdot, x)(u^1\otimes v^1,x)u^2\otimes v^2 = \mathcal{Y}_1(u^1,x)u^2\otimes \mathcal{Y}_2(v^1,x)v^2.$$
\end{theorem}

We now turn our discussions to  the tensor category of $V$-modules. Let $M^1$, $M^2$ be weak $V$-modules. A {\em tensor product} for the ordered pair $(M^1, M^2)$ is a pair $(M, \Y^{12})$ consisting of a weak $V$-module $M$ and an intertwining operator $\Y^{12}$ of type $\left(\begin{tabular}{c}
$M$\\
$M^1$ $M^2$\\
\end{tabular}\right)$ satisfying the following universal property: For any weak $V$-module $W$ and any intertwining operator $I$ of type $\left(\begin{tabular}{c}
$W$\\
$M^1$ $M^2$\\
\end{tabular}\right)$, there exists a unique $V$-homomorphism $\psi$ from $M$ to $W$ such that $I=\psi\circ \Y^{12}$. 

Fix a complex number $z\in \C^{\times}$, for any $\C$-graded $V$-module $M=\oplus_{n\in \C}M_n$ such that $\dim M_n<\infty$ for each $n\in \C$, set $$\overline{M}=\prod_{n\in\C}M_n.$$
For any   $V$-modules $M^1$, $M^2$, $M^3$, a {\em $P(z)$-intertwining
map} $F$ of type $\left(\begin{tabular}{c}
$M^3$\\
$M^1$ $M^2$\\
\end{tabular}\right)$ is a linear map$$F: M^1\otimes M^2\rightarrow \overline{M^3}$$
such that for any $v \in V, w^1\in M^1, w^2\in M^2$,
$$z^{-1}\delta(\frac{x_1-x_0}{z})F(Y_{M^1}(v, x_0)w^1\otimes w^2)+x_0^{-1}\delta(\frac{z-x_1}{-x_0})F(w^1\otimes Y_{M^2}(v, x_1)w^2)$$
$$=x_0^{-1}\delta(\frac{x_1-z}{x_0})Y_{M^3}(v, x_1)F(w^1\otimes w^2).$$

Let $M^1$ and $M^2$ be  $V$-modules.
A {\em $P(z)$-tensor product of $M^1$ and $M^2$} is a pair $$(M^1\boxtimes_{P(z)} M^2, \boxtimes_{P(z)})$$ consisting of a $V$-module $M^1\boxtimes_{P(z)} M^2$  and  a $P(z)$-intertwining
map $\boxtimes_{P(z)}$ of type $\left(\begin{tabular}{c}
$M^1\boxtimes_{P(z)} M^2$\\
$M^1$ $M^2$\\
\end{tabular}\right)$ satisfying the following universal property: For any pair $(M^3,  F)$ consisting of a $V$-module $M^3$ and a $P(z)$-intertwining
map $ F$ of type $\left(\begin{tabular}{c}
$M^3$\\
$M^1$ $M^2$\\
\end{tabular}\right)$, there is a unique $V$-module morphism $\eta$ from $M^1\boxtimes_{P(z)} M^2$ to $M^3$ such that $G=\overline{\eta}\circ F$, where $\overline{\eta}$ is the map from $\overline{M^3}$ to $\overline{M^4}$ uniquely extending $\eta$. The $V$-module $M^1\boxtimes_{P(z)} M^2$ is called a {\em $P(z)$-tensor product module} of $M^1$ and $M^2$. 

We now recall from \cite{HL4} the precise connection between intertwining operators and $P(z)$-intertwining maps of the same type. Fix an integer $p$. For any $V$-modules $M^1, M^2, M^3$, let $\Y$ be an intertwining
operator of type $\left(\begin{tabular}{c}
$M^3$\\
$M^1$ $M^2$\\
\end{tabular}\right)$. We have a linear map $F_{\Y, p}:M^1\otimes M^2\to \overline{M^3}$ given by $$F_{\Y, p}(w^1\otimes w^2)=\Y(w^1, e^{l_p(z)})w^2$$ for any $w^1\in M^1, w^2\in M^2$, where $l_p(z)=\log z+2p\pi i.$ Conversely, given a $P(z)$-intertwining
map $F$ of type $\left(\begin{tabular}{c}
$M^3$\\
$M^1$ $M^2$\\
\end{tabular}\right)$, for any homogeneous elements $w^1\in M^1$, $w^2\in M^2$ and complex number $n\in \C$, let $w^1_nw^2$  be the projection of the image of $w^1\otimes w^2$ under $F$ to the homogeneous subspace of $M^3$ of weight $\wt w^1+\wt w^2-n-1$, multiplied by $e^{(n+1)l_p(z)}$. Set $$\Y_{F, p}(w^1, x)w^2=\sum_{n\in \C}w^1_nw^2x^{-n-1},$$ we obtain a linear map $$M^1\otimes M^2\to M^3\{x\}$$$$w^1\otimes w^2\mapsto \Y_{F, p}(w^1, x)w^2.$$
The following theorem was established in \cite{HL3}, \cite{HL4}.
\begin{theorem}\label{iso}
For $p\in \Z$, the correspondence $\Y\mapsto F_{\Y, p}$ is a linear isomorphism from the vector space  of intertwining
operators of type $\left(\begin{tabular}{c}
$M^3$\\
$M^1$ $M^2$\\
\end{tabular}\right)$ to the vector space  of $P(z)$-intertwining
maps of type $\left(\begin{tabular}{c}
$M^3$\\
$M^1$ $M^2$\\
\end{tabular}\right)$. Its inverse map is given by $F\mapsto \Y_{F, p}$.
\end{theorem}
In the following, we take $p=0$ and the meaning of the intertwining operator associated to a $P(z)$-intertwining map is clear.

In general, the existence of a $P(z)$-tensor product is not obvious. If $V$ is a  vertex operator algebra satisfying  conditions (1) and (2), and let $M^1, ..., M^k$ be the complete list of irreducible $V$-modules. It was proved in a series of papers \cite{HL1}, \cite{HL2}, \cite{HL3}£¬ \cite{H1} and \cite{H3} that the $P(z)$-tensor product $(M^1\boxtimes_{P(z)} M^2,\boxtimes_{P(z)})$ exists for any two $V$-modules $M^1$, $M^2$ and $M^1\boxtimes_{P(z)} M^2\cong\oplus_iN_{M^1, M^2}^{M^i}M^i$.  From  Theorem \ref{iso}, there is a unique intertwining operator $\Y_{P(z);M^1, M^2}$ of type $\left(\begin{tabular}{c}
$M^1\boxtimes_{P(z)} M^2$\\
$M^1$ $M^2$\\
\end{tabular}\right)$ associated to the $P(z)$-intertwining map $\boxtimes_{P(z)}$. From the definition of $P(z)$-tensor product, it is clear that $(M^1\boxtimes_{P(z)} M^2, \Y_{P(z);M^1, M^2})$ is a tensor product of $M^1$ and $M^2$. In the following,  we denote the $P(1)$-tensor product $(M^1\boxtimes_{P(1)} M^2, \boxtimes_{P(1)})$ by $(M^1\boxtimes M^2, \boxtimes_{P(1)})$ and $\Y_{P(1);M^1, M^2}$ by $\Y_{M^1, M^2}.$ For any $w^1\in M^1$ and $w^2\in M^2$, the {\em $P(z)$-tensor product} of $w^1$ and $w^2$ is defined to be $$w^1\boxtimes_{P(z)}w^2=\Y_{P(z);M^1, M^2}(w^1, x)w^2|_{x^n =e^{n\log z}}, n\in \C.$$

 \vskip0.25cm
We now recall some facts about the { braiding isomorphism}  from \cite{H6}. Let $z_1, z_2\in \C^{\times}$ and $\gamma$ be a path in $\C^{\times}$ from $z_1$ to $z_2$. The {\em parallel isomorphism } $\mathcal{T}_{\gamma}:M^1\boxtimes_{P(z_1)}M^2\to M^1\boxtimes_{P(z_2)}M^2$ (see \cite{HL4}, \cite{H6}) is defined as follows: Let $\Y_{P(z_2);M^1, M^2}$ be the intertwining operator associated to the $P(z_2)$-tensor product $(M^1\boxtimes_{P(z_2)} M^2,\boxtimes_{P(z_2)})$ and $l(z_1)$ be the value of logarithm of $z_1$ determined uniquely by $\log z_2$ and the path $\gamma$. Then $\mathcal{T}_{\gamma}$ is characterized by $$\overline{\mathcal{T}_{\gamma}}(w^1\boxtimes_{P(z_1)}w^2)=\Y_{P(z_2);M^1, M^2}(w^1, x)w^2|_{x^n=e^{nl(z_1)}}, n\in \C$$ for any $w^1\in M^1$ and $w^2\in M^2$.
 The parallel isomorphism depends only on the homotopy class of $\gamma$.

For any $V$-modules $M^1, M^2$, the {\em braiding isomorphism} (also called the {\em  commutativity isomorphism}) $${C}_{P(z);M^1, M^2}:M^1\boxtimes_{P(z)} M^2\to M^2\boxtimes_{P(z)} M^1$$ is characterized as follows: For any $w^1\in M^1$ and $w^2\in M^2$, $$\overline{{C}_{P(z);M^1, M^2}}(w^1\boxtimes_{P(z)}w^2)=e^{zL(-1)}\overline{\mathcal{T}_{\gamma_z^-}}(w^2\boxtimes_{P(-z)}w^1),$$where $\gamma_z^-$ is a path from $-z$ to $z$ in the closed upper half plane with $0$ deleted and $\mathcal{T}_{\gamma_z^-}$ is the corresponding parallel isomorphism. When $z=1$, denote the braiding isomorphism ${C}_{P(1);M^1, M^2}:M^1\boxtimes M^2\to M^2\boxtimes M^1$ by ${C_{M^1, M^2}}$.

Let $\Y_{M^1, M^2}$, $\Y_{M^2, M^1}$ be the intertwining operators associated to the $P(1)$-tensor product $(M^1\boxtimes M^2, \boxtimes_{P(1)})$, $(M^2\boxtimes M^1, \boxtimes_{P(1)})$, respectively. Then $(M^1\boxtimes M^2, \Y_{M^1, M^2})$ is  a tensor product of $M^1$ and $M^2$, $(M^2\boxtimes M^1, \Y_{M^2, M^1})$ is  a tensor product of $M^2$ and $M^1$. Consider the intertwining  operator $\Y_{M^2, M^1}^*$ of type $\left(\begin{tabular}{c}
$M^2\boxtimes M^1$\\
$M^1$ $M^2$\\
\end{tabular}\right)$ which is defined by $$\Y_{M^2, M^1}^*(w^1, x)w^2=e^{xL(-1)}\Y_{M^2, M^1}(w^2, -x)w^1.$$
\begin{lemma}\label{eqcomm}
For any $V$-modules $M^1$, $M^2$ and $w^1\in M^1$, $w^2\in M^2$, we have
$$\overline{C_{M^1, M^2}}(\Y_{M^1, M^2}(w^1, x)w^2|_{x^n=e^{n\log 1}})=\Y_{M^2, M^1}^*(w^1, -x)w^2|_{x^n=e^{n\pi i}}.$$
\end{lemma}
 \pf By direct calculation, for any $w^1\in M^1$ and $w^2\in M^2$,
 \begin{align*}
 &\overline{C_{M^1, M^2}}(\Y_{M^1, M^2}(w^1, x)w^2|_{x^n=e^{n\log 1}})\\
 &\ \ =\overline{C_{M^1, M^2}}(w^1\boxtimes w^2)\\
 &\ \ =e^{L(-1)}\overline{\mathcal{T}_{\gamma_1^-}}(w^2\boxtimes_{P(-1)} w^1)\\
 &\ \ =e^{-xL(-1)}\Y_{M^2, M^1}(w^2, x)w^1|_{x^n=e^{n\pi i}}\\
 &\ \ =\Y_{M^2, M^1}^*(w^1, -x)w^2|_{x^n=e^{n\pi i}},
 \end{align*}
 as desired.
 \qed

For any $V$-modules $M^1, M^2, M^3$, the associativity isomorphism $$\ca_{M^1, M^2, M^3}: M^1\bt( M^2\bt M^3)\to (M^1\bt M^2)\bt M^3,$$ was  defined in \cite{H1}, \cite{H2}. Together with the braiding isomorphism defined above, the category $\cc_V$ of $V$-module was proved to be a modular tensor category (see \cite{H6}).
\begin{theorem}\label{kvoa1}
 Let $V$ be a  vertex operator algebra satisfying conditions $(1)$ and $(2)$. Then the $V$-module category $(\cc_V, \boxtimes)$ is a modular tensor category. For $W\in \cc$,  the dual of $W$ is isomorphic to $W'$ and  $\theta_{W}=e^{2\pi iL(0)}$.
 \end{theorem}
Recall that a vertex operator algebra $U$ is called an {\em extension vertex operator algebra} of $V$ if $V$ is a vertex operator subalgebra of $U$ such that $V$, $U$ have the same conformal vector. The following theorem was established in \cite{HKL}.
\begin{theorem}\label{kvoa2} Let $V$ be a  vertex operator algebra satisfying conditions $(1)$ and $(2)$. If $U$ is an extension vertex operator algebra of $V$, then
$U$ induces a commutative algebra $A_{U}$ in $\cc_V$ such that $A_U$ is isomorphic to $ U$ as $V$-module. Conversely, if $U$ is a $V$-module having integral conformal weight and $U$ is a commutative algebra  in $\cc_V$, then $U$ has a vertex operator algebra structure such that $U$ is an extension vertex operator algebra of $V$. 
\end{theorem}
The following result follows immediately from Theorems \ref{kc1}, \ref{semisimple}, \ref{kvoa1}.
\begin{theorem}\label{kvoa5} Let $V$ be a  vertex operator algebra satisfying  conditions $(1)$, $(2)$, 
 $U$ be an extension vertex operator algebra of $V$ such that $U$ is simple and that $V$ appears in $U$ with multiplicity $1$. Then $U$ induces a connected etale algebra $A_{U}$ in $\cc_V$ such that $A_U$ is isomorphic to $ U$ as $V$-module.
\end{theorem}
\subsection{Tensor category of tensor product VOA} Let $V_1$, $V_2$ be   vertex operator algebras satisfying conditions (1) and (2), $(\cc_{V_1}, \bt^1)$, $(\cc_{V_2}, \bt^2)$ be the module categories of $V_1$, $V_2$, respectively. From Theorem \ref{kvoa3},  $V^1\otimes V^2$ satisfies  conditions (1) and $(2)$, then the $V^1\otimes V^2$-module category $(\cc_{V^1\otimes V^2}, \bt^{12})$ is a modular tensor category. Furthermore, we have
\begin{lemma}\label{kvoa7}
For any $V^1$-modules $M^1, M^2$ and $V^2$-modules $N^1, N^2$, the tensor product $(M^1\ot N^1)\bt^{12}(M^2\ot N^2)$ of $M^1\ot N^1$ and $M^2\ot N^2$ is isomorphic to $(M^1\bt^1M^2)\ot (N^1\bt^2 N^2)$. Moreover, under this isomorphism the braiding  isomorphism $$C_{M^1\ot N^1,M^2\ot N^2}:(M^1\ot N^1)\bt^{12}(M^2\ot N^2)\to (M^2\ot N^2)\bt^{12}(M^1\ot N^1)$$ is equal to $C_{M^1, M^2}\ot C_{N^1, N^2}$.
\end{lemma}
\pf Let $\Y_{M^1, M^2}$, $\Y_{N^1, N^2}$  be the intertwining operators associated to the $P(1)$-tensor product $(M^1\boxtimes^1 M^2, \boxtimes^1_{P(1)})$, $(N^1\boxtimes^2 N^2, \boxtimes^2_{P(1)})$, respectively. Let $\sigma(\Y_{M^1, M^2}\ot \Y_{N^1, N^2})$ be
the
image of $\Y_{M^1, M^2}\ot \Y_{N^1, N^2}$ under the isomorphism $\sigma$ defined in Theorem \ref{kvoa6}, then $\sigma(\Y_{M^1, M^2}\ot \Y_{N^1, N^2})$ is an intertwining
operators of type $$\left(\begin{tabular}{c}
$(M^1\bt^1 M^2)\ot (N^1\bt^2 N^2)$\\
$M^1\ot N^1$ $M^2\ot N^2$\\
\end{tabular}\right)$$
From the universal property of $\Y_{M^1, M^2}$, $\Y_{N^1, N^2}$  and  Theorems \ref{kvoa3},  \ref{kvoa6}, it is clear  that $\sigma(\Y_{M^1, M^2}\ot \Y_{N^1, N^2})$ satisfies the universal property in the definition of tensor product, implies that $((M^1\bt^1 M^2)\ot (N^1\bt^2 N^2), \sigma(\Y_{M^1, M^2}\ot \Y_{N^1, N^2}))$ is a tensor product of $M^1\ot N^1$ and $M^2\ot N^2$, then there is a $V^1\ot V^2$-module isomorphism $$\phi_{M^1, M^2;N^1, N^2}:(M^1\bt^1M^2)\ot (N^1\bt^2 N^2)\to (M^1\ot N^1)\bt^{12}(M^2\ot N^2)
 $$such that $$\mathcal{Y}_{M^1\ot N^1,M^2\ot N^2}=\phi_{M^1, M^2;N^1, N^2}\circ \sigma(\Y_{M^1, M^2}\ot \Y_{N^1, N^2}).$$

 To prove the second part of the lemma, note that the conformal vector $\w$ of $V^1\ot V^2$ is equal to $\w_1\ot \1+\1\ot \w_2$ and $L(-1)=\w(0)=\w_1(0)\ot 1+1\ot \w_2(0)=L_1(-1)\ot 1+1\ot L_2(-1)$. Therefore, $e^{xL(-1)}=e^{x(L_1(-1)\ot 1)}e^{x(1\ot L_2(-1))}$.
Then, by Lemma \ref{eqcomm}, we have for any $u^1\in M^1,u^2\in M^2$ and $w^1\in N^1,w^2\in N^2$,
\begin{align*}
&\overline{\phi_{M^2, M^1;N^2, N^1}\circ C_{M^1, M^2}\ot C_{N^1, N^2}\circ\phi_{M^1, M^2;N^1, N^2}^{-1}}((u^1\ot w^1)\boxtimes^{12} (u^2\ot w^2))\\
&\ \ =\overline{\phi_{M^2, M^1;N^2, N^1}\circ C_{M^1, M^2}\ot C_{N^1, N^2}\circ\phi_{M^1, M^2;N^1, N^2}^{-1}}\\
&\ \ \ \ \ \ \ \ \ \ \ \ \ \ \ \ \ \ \ \ \ \ \ \ \ \ \ \   (\mathcal{Y}_{M^1\ot N^1,M^2\ot N^2}(u^1\ot w^1, x)u^2\ot w^2|_{x^n=e^{n\log 1}})\\
&\ \ =\overline{\phi_{M^2, M^1;N^2, N^1}\circ C_{M^1, M^2}\ot C_{N^1, N^2}}\\
&\ \ \ \ \ \ \ \ \ \ \ \ \ \ \ \ \ \ \ \ \ \ \ \ \ \ \ \  (\phi_{M^1, M^2;N^1, N^2}^{-1}\circ\mathcal{Y}_{M^1\ot N^1,M^2\ot N^2}(u^1\ot w^1, x)u^2\ot w^2|_{x^n=e^{n\log 1}})\\
&\ \ =\overline{\phi_{M^2, M^1;N^2, N^1}\circ C_{M^1, M^2}\ot C_{N^1, N^2}}\\
&\ \ \ \ \ \ \ \ \ \ \ \ \ \ \ \ \ \ \ \ \ \ \ \ \ \ \ \  (\sigma(\Y_{M^1, M^2}\ot \Y_{N^1, N^2})(u^1\ot w^1, x)u^2\ot w^2|_{x^n=e^{n\log 1}})\\
&\ \ =\overline{\phi_{M^2, M^1;N^2, N^1}\circ C_{M^1, M^2}\ot C_{N^1, N^2}}\\
&\ \ \ \ \ \ \ \ \ \ \ \ \ \ \ \ \ \ \ \ \ \ \ \ \ \ \ \  ((\Y_{M^1, M^2}(u^1,x)u^2)\ot (\Y_{N^1, N^2}( w^1, x) w^2)|_{x^n=e^{n\log 1}})\\
&\ \ =\overline{\phi_{M^2, M^1;N^2, N^1}}((\Y^*_{M^2, M^1}(u^1,-x)u^2)\ot (\Y^*_{N^2, N^1}( w^1, -x) w^2)|_{x^n=e^{n\pi i}})\\
&\ \ =\overline{\phi_{M^2, M^1;N^2, N^1}}(\sigma(\Y_{M^2, M^1}^*\ot \Y_{N^2, N^1}^*)(u^1\ot w^1, -x)u^2\ot w^2|_{x^n=e^{n\pi i}})\\
&\ \ =\mathcal{Y}_{M^2\ot N^2,M^1\ot N^1}^*(u^1\ot w^1, -x)u^2\ot w^2|_{x^n=e^{n\pi i}}\\
&\ \ =\overline{C_{M^1\ot N^1,M^2\ot N^2}}((u^1\ot w^1)\boxtimes^{12} (u^2\ot w^2)).
\end{align*}
Thus
$$C_{M^1\ot N^1,M^2\ot N^2}=\phi_{M^2, M^1;N^2, N^1}\circ C_{M^1, M^2}\ot C_{N^1, N^2}\circ\phi_{M^1, M^2;N^1, N^2}^{-1}.$$
This finishes the proof.
\qed
\begin{lemma}\label{kvoa8}
For any $V^1$-module $M$ and $V^2$-module $N$, the dual $(M\ot N)'$ of $M\ot N$ is isomorphic to $M'\ot N'$.
\end{lemma}
\pf This follows immediately from Lemma \ref{kvoa7} and the uniqueness of the dual of $M\ot N$.\qed
\begin{lemma}\label{kvoa9}
For any $V^1$-module $M$ and $V^2$-module $N$, we have $$\FPdim_{\cc_{V^1\ot V^2}} M\ot N=\FPdim_{\cc_{V^1}} M \cdot \FPdim_{\cc_{ V^2}} N.$$
\end{lemma}
\pf For any $V^1$-module $M$ and $V^2$-module $N$, define $\phi(M\ot N)=\FPdim_{\cc_{V^1}} M \cdot \FPdim_{\cc_{ V^2}} N$. From Lemma \ref{kvoa7}, $\phi$ is a homomorphism from the Grothendieck ring of $\cc_{V^1\ot V^2}$ to real numbers sending
each isomorphism class to a nonnegative real number. Since homomorphism having this property is unique, we have $\FPdim_{\cc_{V^1\ot V^2}} M\ot N=\FPdim_{\cc_{V^1}} M \cdot \FPdim_{\cc_{ V^2}} N.$ \qed
 \section{Mirror extension of vertex operator algebras}
\def\theequation{3.\arabic{equation}}
\setcounter{equation}{0}
In this section, we shall prove our main results. Let $V^1$, $V^2$ be  vertex operator algebras satisfying conditions (1) and (2). Denote  the module categories of $V^1$, $V^2$ by  $\cc_{V^1}$, $\cc_{V^2}$, respectively. If $V$ is an extension vertex operator algebra of $V^1\otimes V^2$, it is known that $V^1\ot V^2$ is rational by  Theorem \ref{kvoa3},  then  $V$ has the following decomposition as $V^1\ot V^2$-module:
$$V=\oplus_{i\in I, j\in J}Z_{i,j}M^i\otimes N^j,$$
where $Z_{i, j}$, $(i\in I, j\in J)$,  are nonnegative integers  and $\{M^i|i\in I\}$ (resp. $\{N^j|j\in J\}$) are inequivalent irreducible $V^1$-modules (resp. $V^2$-modules) such that $Z_{i, j}\neq 0$  for some $(i,j)\in I\times J$. 
In the following we  assume that $V$ is an  extension vertex operator algebra  of $V^1\ot V^2$ such that $V$ is simple, self dual
and $\Hom_{V^1\otimes V^2}(V^1\otimes N^j, V)=\C$ (resp. $\Hom_{V^1\otimes V^2}(M^i\otimes V^2, V)=\C$) if and only if $N^j=V^2$ (resp. $M^i=V^1$).
Then we have

\begin{theorem}\label{c1}
(i) $Z_{i,j}=1$ if $Z_{i,j}\neq 0$. Moreover, this induces a bijection $\tau$ from $I$ to $J$ such that $\tau(i)=j$ if and only if $Z_{i, j}=1$. \\
(ii) The set $\{M^i| i\in I\}$ $($resp. $\{N^j| j\in J\})$ is closed under the tensor product of category $\cc_{V^1}$ $($resp. $\cc_{V^2})$. \\
(iii) For any $i_1, i_2,i_3\in I$,  $N_{M^{i_1}, M^{i_2}}^{M^{i_3}}=N_{(N^{\tau(i_1)})', (N^{\tau(i_2)})'}^{(N^{\tau(i_3)})'}$.
\end{theorem}
\pf Let $\{M^i| i\in I^c\}$ $($resp. $\{N^j| j\in J^c\})$ be irreducible $V^1$-modules (resp. $V^2$-modules) such that $\{M^i| i\in I\}$ $($resp. $\{N^j| j\in J\})$ is a subset of $\{M^i| i\in I^c\}$ $($resp. $\{N^j| j\in J^c\})$ and is closed under the tensor product of category $\cc_{V^1}$ $($resp. $\cc_{V^2})$. For $i\in I^c, j\in J^c$, we take $Z_{i, j}=0$ if $i\in I^c\backslash I$ or $j\in J^c\backslash J$, then we have $V=\oplus_{i\in I^c, j\in J^c}Z_{i,j}M^i\otimes N^j.$

Consider the module category $\cc_{V^1\otimes V^2}$ of $V^1\ot V^2$,
 from Theorem \ref{kvoa5}, $V$ induces an etale algebra $A$ in $\cc_{V^1\otimes V^2}$ such that $A$ is isomorphic to $ V$ as $V^1\otimes V^2$-module. From Proposition \ref{kc3}, Theorem \ref{kc2} and Lemmas \ref{kvoa7}, \ref{kvoa8}, we have for any $i_1, i_2\in I^c$,
\begin{align*}
&\Hom_A(A\bt^{12} (M^{i_1}\ot V^2), A\bt^{12} (M^{i_2}\ot V^2))\\
&\ \ =\Hom(M^{i_1}\ot V^2, A\bt^{12}(M^{i_2}\ot V^2))\\
&\ \ =\Hom((M^{i_1}\ot V^2)\bt^{12} (M^{i_2}\ot V^2)',  A)\\
&\ \ =\Hom(( M^{i_1}\bt^1 (M^{i_2})')\ot V^2,  A),
\end{align*}
it follows from the assumption that $\Hom_A(A\bt^{12}(M^{i_1}\ot V^2), A\bt^{12}(M^{i_2}\ot V^2))=\C$ if and only if $M^{i_1}=M^{i_2}.$ This implies that $\{A\bt^{12}(M^{i}\ot V^2)|i\in I^c\}$ are simple left $A$-modules and mutually inequivalent. Similarly, we have $\Hom_A(A\bt^{12}(V^1\ot N^{j_1}), A\bt^{12}(V^1\ot N^{j_2}))=\C$ if and only if $N^{j_1}=N^{j_2}.$ Therefore, $\{A\bt^{12}(V^1\ot N^{j})|j\in J^c\}$ are simple left $A$-modules and mutually inequivalent.

Also from Proposition \ref{kc3}, Theorem \ref{kc2} and Lemmas \ref{kvoa7}, \ref{kvoa8}, we have for any $i\in I^c$ and $j\in J^c$,
\begin{align*}
&\Hom_A(A\bt^{12}(M^{i}\ot V^2),  A\bt^{12}(V^1\ot N^{j}))\\
&\ \ =\Hom((M^{i}\ot V^2), A\bt^{12}( V^1\ot N^{j}))\\
&\ \ =\Hom((M^{i}\ot V^2)\bt^{12} (V^1\ot N^{j})',  A)\\
&\ \ =\Hom(M^{i}\ot (N^{j})',  A),
\end{align*}
implies that 
$\dim \Hom_A(A\bt^{12}(M^{i}\ot V^2), A\bt^{12}( V^1\ot (N^{j})'))=Z_{i, j}$ and that  $Z_{i, j}=1$ if  $Z_{i, j}\neq 0$. Similarly, we have $\dim \Hom_A(A\bt^{12}( V^1\ot (N^{j})'), A\bt^{12}(M^{i}\ot V^2))=Z_{i', j'}$ and $Z_{i', j'}=1$ if $Z_{i', j'}\neq 0$, where we use $N^{j'}$ to denote $(N^j)'$. Note that $V$ is self dual, therefore $Z_{i, j}=Z_{i', j'}$. This implies that  $A\bt^{12}(M^{i}\ot V^2)$ is isomorphic to $A\bt^{12}( V^1\ot (N^{j})')$ as left $A$-module if $Z_{i, j}\neq 0$ and that  $A\bt^{12}(M^{i}\ot V^2)$ is not isomorphic to any one of $\{A\bt^{12}( V^1\ot N^{j})|j\in J^c\}$ if $i\in I^c\backslash I$. Since $\{A\bt^{12}(M^{i}\ot V^2)| i\in I^c\}$ (resp. $\{A\bt^{12}( V^1\ot N^{j})|j\in J^c\}$) are inequivalent left $A$-modules, for any $i\in I$ there is a unique $j\in J$ such that $Z_{i, j}\neq 0$. Therefore, we have a bijection $\tau$ from $I$ to $J$ such that $\tau(i)=j$ if and only if $Z_{i, j}\neq 0$.

For any $i\in I, j\in J$, consider the left $A$-module  $A\bt^{12}(M^i\ot (N^j)')$, from the discussion above and Lemmas \ref{kvoa7}, \ref{kvoa8}, we have
\begin{align*}
\begin{split}
A\bt^{12}(M^i\ot (N^j)')&\cong A\bt^{12}((M^i\ot V^2)\bt^{12} (V^1\ot (N^j)' ))\\
&\cong (A\bt^{12} (M^i\ot V^2))\bt^{12} (V^1\ot (N^j)' )\\
&\cong(A\bt^{12}(V^1\ot (N^{\tau(i)})'))\bt^{12}(V^1\ot (N^j)')\\
&\cong A\bt^{12}((V^1\ot (N^{\tau(i)})')\bt^{12} (V^1\ot (N^j)'))\\
&\cong A\bt^{12}(V^1\ot ((N^{\tau(i)})'\bt^2 (N^j)'))\\
&\cong\oplus_lN_{(N^{\tau(i)})',(N^j)'}^{N^l}A\bt^{12}(V^1\ot N^l).
\end{split}
\end{align*}
On the other hand,
\begin{align*}
\begin{split}
A\bt^{12}(M^i\ot (N^j)')&\cong A\bt^{12}((V^1\ot (N^j)')\bt^{12} (M^i\ot V^2))\\
&\cong (A\bt^{12}(V^1\ot (N^j)'))\bt^{12} (M^i\ot V^2)\\
&\cong  (A\bt^{12}(M^{\tau^{-1}(j)}\ot V^2 ))\bt^{12}(M^i\ot V^2)\\
&\cong  A\bt^{12}((M^{\tau^{-1}(j)}\ot V^2 )\bt^{12}(M^i\ot V^2))\\
&\cong A\bt^{12}((M^{\tau^{-1}(j)}\bt^1 M^i)\ot V^2 )\\
&\cong \oplus_{k}N_{M^{\tau^{-1}(j)}, M^i}^{M^k}A\bt^{12}(M^k\ot V^2 ).
\end{split}
\end{align*}
Comparing the two decompositions, we obtain that the set $\{M^i|i\in I\}$ (resp. $\{N^j|j\in J\}$) is closed under the tensor product of category $\cc_{V_1}$ (resp. $\cc_{V_2}$)  and that $N_{M^{\tau^{-1}(j)}, M^i}^{M^{\tau^{-1}(l')}}=N_{(N^j)', (N^{\tau(i)})'}^{N^l}$. Note that $\tau(i')=\tau(i)'$, then we have
$N_{M^{i_1}, M^{i_2}}^{M^{i_3}}=N_{(N^{\tau(i_1)})', (N^{\tau(i_2)})'}^{(N^{\tau(i_3)})'}=N_{N^{\tau(i_1)}, N^{\tau(i_2)}}^{N^{\tau(i_3)}}$.  The proof is complete.
\qed
\vskip0.5cm
From Theorem \ref{c1}, we have $\{M^i|i\in I\}$ (resp. $\{N^j|j\in J\}$) is closed under the tensor product of $\cc_{V^1}$ (resp. $\cc_{V^2}$). Let $\cc_{V^1}^0$ (resp. $\cc_{V^2}^0$) be the fusion subcategory of $\cc_{V^1}$ (resp. $\cc_{V^2}$) such that the Grothendieck ring of $\cc_{V^1}^0$ (resp. $\cc_{V^2}^0$) is equal to  $\{M^i|i\in I\}$ (resp. $\{N^j|j\in J\}$).  $\cc_{V^1}^0$ (resp. $\cc_{V^2}^0$) can be viewed as the subcategory of  $\cc_{V^1\otimes V^2}$ by identifying $M^i$ (resp. $N^j$) with $M^i\ot V^2$ (resp. $V^1\ot N^j$). Note that $\{M^i\ot N^j|i\in I, j\in J\}$ is closed under the tensor product of $\cc_{V^1\otimes V^2}$. Let $\cc$ be the fusion subcategory of $\cc_{V^1\ot V^2}$ such that the Grothendieck ring of $\cc$ is equal to $\{M^i\ot N^j|i\in I, j\in J\}$. Therefore, 
$A$ is a connected etale algebra in $\cc$. Let $\cc_{A}$ be the category of $A$-module in $\cc$. Then we have
\begin{lemma}\label{c2}
$\FPdim \cc_A=\FPdim \cc_{V^1}^0=\FPdim \cc_{V^2}^0$.
\end{lemma}
\pf
From Theorem \ref{c1}, we get for any $i\in I$, $\FPdim_{\cc_{V^1}^0}M^i=\FPdim_{\cc_{V^2}^0}(N^{\tau(i)})'=\FPdim_{\cc_{V^2}^0}N^{\tau(i)}$, thus $\FPdim{\cc_{V^1}^0}=\FPdim{\cc_{V^2}^0}$. It follows from Lemma \ref{kvoa9} that $\FPdim_{\cc}A=\FPdim \cc_{V^1}^0=\FPdim \cc_{V^2}^0 $. Note that $\FPdim \cc=\FPdim \cc_{V^1}^0\cdot \FPdim \cc_{V^2}^0$, from Proposition \ref{kc4}  we have $\FPdim \cc_A=\FPdim \cc_{V^1}^0=\FPdim \cc_{V^2}^0$.
\qed

The following result plays an  important role in the proof of our main results.
\begin{theorem}\label{kvoa11}
There is a braid-reversing equivalence $$\ct: \cc_{V^1}^0\to \cc_{V^2}^0 $$ such that $\ct$ sends  $M^i\in \cc_{V^1}^0$ to $(N^{\tau(i)})'\in \cc_{V^2}^0.$
\end{theorem}
\pf The proof is similar to that of Theorem 5.1 in \cite{OS}. Let $\cc$ and $\cc_A$ the categories defined above. Consider the free module functor $$\cf_A: \cc\to \cc_A,\ X\mapsto X\bt^{12} A.$$ We now prove that the restriction of $\cf_A$ to the subcategory $\cc_{V^1}^0$ of $\cc$ is fully faithful, that is $\Hom_{\cc_{V^1}^0}(M^1, M^2)=\Hom_{\cc_A}(\cf_A(M^1), \cf_A(M^2))$ for any $M^1, M^2\in \cc_{V^1}^0$. In fact, from Proposition \ref{kc3}, Theorem \ref{kc2} and Lemmas \ref{kvoa7}, \ref{kvoa8}, we have
\begin{align*}
\Hom_{\cc_A}(\cf_A(M^1), \cf_A(M^2))&=\Hom_{\cc_A}((M^1\ot V^2)\bt^{12}A, (M^2\ot V^2)\bt^{12}A)\\
&=\Hom(M^1\ot V^2, (M^2\ot V^2)\bt^{12}A)\\
&=\Hom((M^2\ot V^2)'\bt^{12}(M^1\ot V^2),  A)\\
&=\Hom(((M^2)'\bt^1 M^1)\ot V^2,  A)\\
&=\Hom(((M^2)'\bt^1 M^1)\ot V^2,  V^1\ot V^2)\\
&=\Hom(M^1\ot V^2,  M^2\ot V^2)\\
&=\Hom(M^1,  M^2).
\end{align*}
Similarly, the restriction of $\cf_A$ to the subcategory $\cc_{V^2}^0$ is fully faithful. It follows from Lemma \ref{c2} that $\cf_A$ restricted to $\cc_{V^2}^0$ and $\cc_{V^2}^0$ are monoidal equivalences, and then $\cc_{V^1}^0$ and $\cc_{V^2}^0$ are monoidal equivalent.

Recall that the free module functor $\cf_A$  has the structure of a central functor
given by $X\mapsto  (X \boxtimes^{12} A, C_{-,X})$, where $C_{-,X}$ denotes the braiding isomorphism in $\cc$, then we have two fully faithful braided monoidal functors $\cf_1: \cc_{V^1}^0\to \cz(\cc_A)$ and  $\cf_2: \cc_{V^2}^0\to \cz(\cc_A)$. Furthermore, from Theorem \ref{kvoa10} and Lemma \ref{kvoa7} we have$$\cf_1( \cc_{V^1}^0)\subset (\cf_2( \cc_{V^2}^0))',$$
where $(\cf_2( \cc_{V^2}^0))'$ denotes the centralizer of $\cf_2( \cc_{V^2}^0)$.

On the other hand, since $\cc_A$ is monoidal equivalent to $\cc_{V^2}^0$, we could consider $\cc_A$ as a braided monoidal category with the braiding induced from $\cc_{V^2}^0$. It follows from Theorem \ref{kc5} that there exist fully faithful braided monoidal  functor $i:\cc_A\to \cz(\cc_A)$ and  braid-reversing  monoidal functor $j: \cc_A\to \cz(\cc_A)$ such that $(i(\cc_A))'=j(\cc_A)$. Note that  the image of $i$ is the same as the image of $\cf_2$, hence $(\cf_2(\cc_{V^2}^0))'=j(\cc_A)$. Since $j$ and $\cf_1$ are fully faithful, we have $\FPdim (\cf_2(\cc_{V^2}^0))'=\FPdim j(\cc_A)=\FPdim \cc_A=\FPdim \cc_{V^1}^0= \FPdim \cf_1( \cc_{V^1}^0)$, implies $\cf_1( \cc_{V^1}^0)= (\cf_2( \cc_{V^2}^0))'=j(\cc_A)$. Let $\ct$ be  the following composition functor $$\ct:\cc_{V^1}^0\stackrel{\cf_1}{\to}\cf_1(\cc_{V^1}^0)\stackrel{j^{-1}}{\to} \cc_A \stackrel{i}{\to} \cf_2(\cc_{V^2}^0)\stackrel{\cf_2^{-1}}{\to} \cc_{V^2}^0,$$ where all functors present except for $j^{-1}$ are braided with $j^{-1}$ are braid-reversing. Then we have $\ct:\cc_{V^1}^0\to \cc_{V^2}^0$ is a braid-reversing equivalence.

From Proposition \ref{kc3}, Theorem \ref{kc2} and Lemmas \ref{kvoa7}, \ref{kvoa8}, we have
\begin{align*}
\Hom_{\cc_A}((M^i\ot V^2)\bt^{12} A, (V^1\ot N^j)\bt^{12} A)&=\Hom(M^i\ot V^2, (V^1\ot N^j)\bt^{12} A)\\
&=\Hom((V^1\ot N^j)'\bt^{12} (M^i\ot V^2), A)\\
&=\Hom(M^i\ot (N^j)', A).
\end{align*}
It follows that $\Hom_{\cc_A}((M^i\ot V^2)\bt^{12} A, (V^1\ot N^j)\bt^{12} A)=\C$ if and only if $j=\tau(i)'$. Therefore, $\ct(M^i)=(N^{\tau(i)})'.$ The proof is complete.
\qed

\begin{remark}
When $V^1$ is the affine vertex operator algebra $L_{sl_m}(n,0)$ and $V^2$ is the affine vertex operator algebra $L_{sl_n}(m,0)$, the result in Theorem \ref{kvoa11} has been obtained in \cite{OS}, and  gives an explanation of  the  level-rank duality \cite{F}. From this point of view,  the vertex operator algebra $V^1$ is dual to $V^2$.
\end{remark}
 We are now in a position to prove our main result in this paper.
\begin{theorem}\label{kvoa13}
Let $M^0=V^1, M^1, ..., M^k$ be simple objects in $\cc_{V^1}^0$ such that $$U=M^0\oplus M^1\oplus\cdots \oplus M^k$$ is an extension vertex operator algebra of $V^1$. Then $$U^c=V^2\oplus (N^{\tau(1)})'\oplus\cdots \oplus (N^{\tau(k)})'$$ has a vertex operator algebra structure such that $U^c$ is an extension vertex operator algebra of $V^2$. Moreover, $U^c$ is a simple vertex operator algebra if $U$ is a simple vertex operator algebra.
\end{theorem}
\pf 
It follows immediately from Theorems \ref{kvoa2}, \ref{kvoa11} that $U^c$ induces an algebra $(A_{U^c},\mu)$ in the  category $\cc_{V^2}^0$ such that$$\mu\circ C^{-1}_{U^c,U^c}=\mu.$$
Note that the conformal weight of $U^c$ is integral, by e.q. (2.8) in \cite{HKL}, we have $$C^{-1}_{U^c,U^c}=C_{U^c,U^c}\circ\theta_{U^c\bt^2 U^c},$$implies
\begin{align*}
&\mu=\mu\circ C^{-1}_{U^c,U^c}\\
&\ \ =\mu\circ C_{U^c,U^c}\circ\theta_{U^c\bt^2 U^c}\\
&\ \ =\mu\circ C_{U^c,U^c},
\end{align*}
where in the last equality we have used the facts that $\mu$ is a $V^2$-module morphism and that the conformal weight of $U^c$ is integral. Thus $(A_{U^c},\mu)$ is a commutative algebra in the  category $\cc_{V^2}^0$. Then it follows from Theorem \ref{kvoa2} that there is a vertex operator algebra structure on $U^c$.

Assume that $U$ is a simple vertex operator algebra. We now prove that $U^c$ is a simple vertex operator algebra. Otherwise, $U^c$ has a nontrivial ideal $U^1$, then it follows from Theorem 3.4 in \cite{HKL} that $U^1$ induces a nontrivial $A_{U^c}$-submodule $(U^1,\mu_{U^1})$ of $A_{U^c}$ such that $$\mu_{U^1}\circ C_{U^1, U^c}\circ C_{U^c, U^1}=\mu_{U^1}.$$Thus, by Theorem \ref{kvoa11}, $A_{U}$ has a nontrivial $A_{U}$-submodule $(\tilde U^1, \mu_{\tilde U^1} )$ such that $$\mu_{\tilde U^1}\circ C^{-1}_{\tilde U^1, U}\circ C^{-1}_{U, \tilde U^1}=\mu_{\tilde U^1},$$  where $A_{U}$ is the etale algebra induced from the vertex operator algebra $U$. Since the conformal weight of $U$ is integral, by a similar argument as above, we have
$$\mu_{\tilde U^1}\circ C_{ \tilde U^1, U}\circ C_{U, \tilde U^1}=\mu_{\tilde U^1}.$$It follows from Theorem 3.4 in \cite{HKL} that the vertex operator algebra $U$ has a nontrivial ideal, this is a contradiction. The proof is complete.\qed
\vskip0.5cm
As an application, we  prove the mirror extension conjecture in \cite{DJX}. First, recall that if $(U, Y, \1, \w)$ is a vertex operator algebra and
$(V, Y, \1, \w')$ is a vertex operator subalgebra of $U$ such
that $\w'\in U_2$ and $L(1)\w'=0$, then  $(V^c, Y, \1, \w-\w')$ is a vertex operator subalgebra of
$U$, where $V^c=\{v\in U|L'(-1)v=0\}$ (see \cite{FZ}).
 In the following we assume that $(V, Y, \1, \w)$, $(V^c, Y, \1, \w-\w')$,
$(U, Y, \1, \w')$ satisfy conditions (1), (2)
 and $(V^c)^c=V$, it follows  from Theorem \ref{c1} that $U$ has the decomposition $$U=V\ot V^c\oplus (\oplus_{i=1}^nM^i\ot N^i)$$ such that  $\{M^i|1\leq i\leq n\}$ (resp. $\{N^i|1\leq i\leq n\}$) are mutually inequivalent irreducible $V$-modules (resp. $V^c$-modules). Let $$V^e=V\oplus (\oplus_{i=1}^nm_iM^i)$$ and $$(V^c)^e=V^c\oplus (\oplus_{i=1}^nm_i(N^i)'),$$ where $m_i$, $1\leq i\leq n$, are nonnegative integers. 
 Then we have
\begin{theorem}\label{kvoa12}
 If there is a vertex operator algebra structure on $V^e$ such that $V^e$ is an extension vertex operator algebra of $V$. Then $(V^c)^e$ has a vertex operator algebra structure such that $(V^c)^e$ is an extension vertex operator algebra of $V^c$. Moreover, $(V^c)^e$ is a simple vertex operator algebra if $V^e$ is a simple vertex operator algebra.
\end{theorem}
\pf This follows immediately from Theorem \ref{kvoa13}.\qed

\end{document}